\documentclass[12pt]{amsart}

\usepackage{amscd,amssymb}
\usepackage{amsthm,amsmath,amssymb}

\usepackage[matrix,arrow]{xy}

\newtheorem{theorem}[equation]{Theorem}

\newtheorem{proposition}[equation]{Proposition}
\newtheorem{lemma}[equation]{Lemma}
\newtheorem{corollary}[equation]{Corollary}
\newtheorem{conjecture}[equation]{Conjecture}

\theoremstyle{definition}
\newtheorem{example}[equation]{Example}
\newtheorem{definition}[equation]{Definition}

\theoremstyle{remark}
\newtheorem{remark}[equation]{Remark}

\author{Ivan Cheltsov}

\title{Factorial threefolds and Shokurov vanishing}

\email{cheltsov@yahoo.com}

\thanks{All varieties are assumed to be projective, normal and defined over $\mathbb{C}$.}
\begin{document}

\begin{abstract}
We apply the Shokurov vanishing theorem to prove the factoriality
of the following nodal threefolds: a complete intersection of
hypersurfaces $F$ and $G\subset\mathbb{P}^{5}$ of degree $n$ and
$k$ respectively, where $G$ is smooth, $|\mathrm{Sing}(F\cap
G)|\leqslant(n+k-2)(n-1)/5$, $n\geqslant k$; a double cover of a
smooth hypersurface $F\subset\mathbb{P}^{4}$ of degree $n$
branched over a surface that is cut out on $F$ by a hypersurface
$G$ of degree $2r\geqslant n$, and $|\mathrm{Sing}(F\cap
G)|\leqslant(2r+n-2)r/4$.
\end{abstract}

\maketitle

\section{Introduction.}
\label{section:introduction}

A Weil divisor is a $\mathbb{Q}$-Cartier divisor if some its
nonzero multiple is a Cartier divisor, a variety has
$\mathbb{Q}$-factorial singularities if every Weil divisor on it
is a $\mathbb{Q}$-Cartier divisor, a variety is
$\mathbb{Q}$-fac\-to\-ri\-al if its singularities are
$\mathbb{Q}$-factorial. Smooth varieties are
$\mathbb{Q}$-factorial.

Birational geometry of many singular varieties crucially depend on
the $\mathbb{Q}$-factoriality condition. For example, all
$\mathbb{Q}$-factorial nodal\footnote{A variety is called nodal if
all its singular points are isolated ordinary double points.}
quartic threefolds are known to be nonra\-ti\-onal (see
\cite{IsMa71}, \cite{Pu88b}, \cite{Me03}) and all
$\mathbb{Q}$-factorial double covers of $\mathbb{P}^{3}$ branched
over nodal sextic surfaces are nonrational (see \cite{Is80b},
\cite{Pu97}, \cite{ChPa04}), which is false in the
non-$\mathbb{Q}$-factorial case.

\begin{example}
\label{example:Burkhardt-quartic} Every nodal quartic threefold in
$\mathbb{P}^{4}$ does not have more than $45$ singular point (see
\cite{Va83}, \cite{Fr86}). There are nodal quartic threefolds
having any number of singular points up to $45$, and there is a
unique (see \cite{JSV90}) nodal quartic threefold
$\mathcal{B}_{4}$ with $45$ singular points, which is called a
Burkhardt quartic (see \cite{To36}, \cite{Pet98}) and can be given
as
$$
w^4-w\big(x^3+y^3+z^3+t^3\big)+3xyzt=0\subset \mathbb{P}^4\cong{\mathrm{Proj}}\Big({\mathbb C}[x,y,z,t,w]\Big),%
$$
which implies that $\mathcal{B}_{4}$ is determinantal and
rational. The quartic $\mathcal{B}_{4}$ is a unique invariant
having degree $4$ of the simple group
${\mathrm{PSp}}(4,\mathbb{Z}_{3})$ of order $25920$ (see
\cite{vdGr87}, \cite{Hu96}, \cite{HoWe01}, \cite{HulSa02}), and
singular points of the quartic $\mathcal{B}_{4}$ corresponds to
$45$ tritangents of a smooth cubic surface, which is related to
the fact that the Weil group $E_{6}$ is a nontrivial extension of
the group ${\mathrm{PSp}}(4,\mathbb{Z}_{3})$ by $\mathbb{Z}_{2}$.
The quartic  $\mathcal{B}_{4}$ contains a plane, which is not a
Cartier divisor, because the plane is not cut out on
$\mathcal{B}_{4}$ by any hypersurface. On the other hand, the
local class group of an ordinary double point is $\mathbb{Z}$,
which implies that every nonzero multiple of a plane contained in
$\mathcal{B}_{4}$ is not a Cartier divisor. So, the quartic
$\mathcal{B}_{4}$ is not $\mathbb{Q}$-fac\-to\-ri\-al, one can
show that $\mathrm{Cl}(\mathcal{B}_{4})\cong \mathbb{Z}^{16}$ (see
\cite{HoWe01}), but the Lefschetz implies that
$\mathrm{Pic}(\mathcal{B}_{4})\cong\mathbb{Z}$.
\end{example}

\begin{example}
\label{example:Barth-sextic} Let $\pi:X\to \mathbb{P}^3$ be a
double cover ramified in a Barth sextic surface
$$
4\big(\tau^2x^2-y^2\big)\big(\tau^2y^2-z^2\big)\big(\tau^2z^2-x^2\big)=w^2\big(1+2\tau\big)\big(x^2+y^2+z^2-w^{2}\big)^2\subset\mathbb{P}^{3}\cong\mathrm{Proj}\Big(\mathbb{C}[x,y,z,w]\Big),
$$
where $\tau=(1+\sqrt{5})/2$. Then $X$ is nodal and has $65$
singular points (see \cite{Ba96}), but every nodal sextic surface
has at most $65$ singular points (see \cite{JaRu97}, \cite{Wa98}).
There is a determinantal quartic threefold $V\subset\mathrm{P}^4$
with $42$ ordinary double points such that the diagram
$$
\xymatrix{
&Y\ar@{-->}[d]_{\rho}\ar@{^{(}->}[rr]&&\mathbb{P}^4\ar@{-->}[d]^{\gamma}&\\%
&X\ar@{->}[rr]_{\pi}&&\mathbb{P}^{3}&}
$$
commutes (see  \cite{En99}, \cite{Pet98}), where $\rho$ is a
birational map and $\gamma$ is the projection from an ordinary
double point of the quartic $Y$, which implies the rationality of
$X$, because determinantal quartics are rational. The rational map
$\rho$ can be decomposed as a composition of a blow up of a
singular point of the quartic $Y$ and a consecutive blow down of
the proper transforms of $24$ lines on the quartic $Y$ that pass
through the blown up singular point. Every nonzero multiple of the
image of the exceptional divisor of the blow up of the singular
point of $Y$ is not a Cartier divisor, which implies that $X$ is
not $\mathbb{Q}$-factorial, but one can show that
$\mathrm{Pic}(X)\cong\mathbb{Z}$ and $\mathrm{Cl}(X)\cong
\mathbb{Z}^{14}$ (see Example 3.7 in \cite{En99}).
\end{example}

It is natural to ask how a global topological condition of being
$\mathbb{Q}$-factorial depends on the number of singular points of
a nodal threefold. To illustrate a general picture let us consider
nodal hypersurfaces in $\mathbb{P}^{4}$. Let $V$ be a hypersurface
in $\mathbb{P}^{4}$ of degree $n$ that has at most ordinary double
points. Then the variety $V$ is $\mathbb{Q}$-factorial if and only
if
$$
\mathrm{rk}\,H^{2}\big(V, \mathbb{Z}\big)=\mathrm{rk}\,H_{4}\big(V, \mathbb{Z}\big),%
$$
which is true in the smooth case due to Poincare duality.
Moreover, the following important result holds (see
 \cite{Cl83}, \cite{We87}, \cite{Di90},
\cite{Cy01}).

\begin{proposition}
\label{proposition:defect-II}%
The hypersurface $V$ is $\mathbb{Q}$-factorial if and only if its
singular points impose independent linear conditions on global
sections of the sheaf $\mathcal{O}_{\mathbb{P}^{4}}(2n-5)$.
\end{proposition}

In particular, the hypersurface $V$ is always
$\mathbb{Q}$-factorial when $|\mathrm{Sing}(V)|\leqslant 2n-4$.

\begin{remark}
\label{remark:factoriality-of-hypersurface}%
Let $X$ be either a nodal complete intersection of two
hypersurface in $\mathbb{P}^{5}$, or a nodal double cover of a
smooth hypersurface in $\mathbb{P}^{4}$. Then the group
$\mathrm{Pic}(X)$ is generated either by the class of a hyperplane
section or by a pull back of the class of a hyperplane section.
The threefold $X$ is usually called factorial in the case when a
similar statement holds for the group $\mathrm{Cl}(X)$. However,
the local class group of an isolated ordinary double point is
$\mathbb{Z}$ (see \cite{Mi68}), which implies that the following
conditions are equivalent:
\begin{itemize}
\item the variety $V$ is $\mathbb{Q}$-factorial;%
\item the variety $V$ is factorial;%
\item the isomorphism $\mathrm{Cl}(V)\cong\mathrm{Pic}(V)$ holds;%
\item the isomorphism $\mathrm{Cl}(V)\cong\mathbb{Z}$ holds;%
\item the equality $\mathrm{rk}\,\mathrm{Cl}(V)=1$ holds.
\end{itemize}
\end{remark}

Let us consider the simplest example of the hypersurface $V$ that
is not $\mathbb{Q}$-factorial.

\begin{example}
\label{example:non-Q-factoriality-hypersurface-containing-plane}
Suppose that the hypersurface $V$ is given by the equation
$$
xg\big(x,y,z,t,w\big)+yf\big(x,y,z,t,w\big)=0\subset\mathbb{P}^{4}\cong\mathrm{Proj}\Big(\mathbb{C}[x,y,z,t,w]\Big),
$$
where $g$ and $f$ are general polynomials of degree $n-1$. Then
$V$ is indeed nodal, it contains the plane $x=y=0$ and
$|\mathrm{Sing}(V)|=(n-1)^{2}$, which implies that $V$ is not
$\mathbb{Q}$-factorial.
\end{example}

The problem of $\mathbb{Q}$-factoriality of nodal threefolds is
related to the Shokurov vanishing theorem (see \cite{Sho93},
\cite{Ko91}, \cite{Ko97}, \cite{Am99}). Let us illustrate this
relation on the following example.

\begin{proposition}
\label{proposition:factoriality-answer-to-Ciliberto}%
Let $\mathcal{H}$ be a linear system of hypersurface in
$\mathbb{P}^{4}$ of  degree $k<n/2$ that pass through the points
of the set $\mathrm{Sing}(V)$, and let
$\hat{\mathcal{H}}=\mathcal{H}\vert_{V}$. Suppose that the base
locus of the linear system $\hat{\mathcal{H}}$ is
zero-dimensional. Then the hypersurface $V$ is
$\mathbb{Q}$-factorial.
\end{proposition}

\begin{proof}
Let $P$ ba an arbitrary singular point of $V$. Then it follows
from Proposition~\ref{proposition:defect-II} that in order to
conclude the proof we must find a hypersurface in $\mathbb{P}^{4}$
of degree $2n-5$ that pass through all points of the set
$\mathrm{Sing}(V)\setminus P$ and does not pass through the point
$P$.

Suppose that the base locus of the linear system ${\mathcal{H}}$
is zero-dimensional. Let $\Lambda$ be a base locus of
$\mathcal{H}$. Then $\mathrm{Sing}(V)\subseteq\Lambda$. Take
sufficiently general divisors $H_{1},\ldots, H_{s}$ in the linear
system $\mathcal{H}$ for $s\gg 0$. Put $X=\mathbb{P}^{4}$,
$B_{X}={\frac{4}{s}}\sum_{i=1}^{s}H_{i}$ and
$$
\mathrm{Sing}\big(V\big)\setminus P=\big\{P_{1},\ldots, P_{r}\big\},%
$$
where $P_{i}$ is a point. Let $f:V\to X$ be a blow up of all
points in $\mathrm{Sing}(V)\setminus P$. Then
$$
K_{V}+B_{V}+\sum_{i=1}^{r}\Big(\mathrm{mult}_{P_{i}}\big(B_{X}\big)-4\Big)E_{i}+f^{*}\big(H\big)=f^{*}\Big(\big(4k-4\big)H\Big)-\sum_{i=1}^{r}E_{i},%
$$
where $E_{i}$ is the $f$-exceptional divisor such that
$f(E_{i})=P_{i}$, the divisor $B_{V}$ is a proper transform of the
the divisor $B_{X}$ on the variety $V$, and $H$ is a hyperplane in
$\mathbb{P}^{4}$. Let
$$
\hat{B}_{V}=B_{V}+\sum_{i=1}^{r}\Big(\mathrm{mult}_{P_{i}}\big(B_{X}\big)-4\Big)E_{i},
$$
and $\bar{P}$ be a point of the variety $V$ such that
$f(\bar{P})=(P)$. Then the divisor $\hat{B}_{V}$ is effective,
because $\mathrm{mult}_{P_{i}}(B_{X})\geqslant 4$ for every $i$.
We have $\mathrm{mult}_{\bar{P}}(B_{V})\geqslant 4$, which implies
that $\bar{P}$ is an isolated center of log canonical
singularities of the log pair $(V, \hat{B}_{V})$. Hence, the map
$$
H^{0}\Big(\mathcal{O}_{V}\Big(f^{*}\big((4k-4)H\big)-\sum_{i=1}^{r}E_{i}\Big)\Big)\to H^{0}\Big(\mathcal{O}_{\mathcal{L}\big(V,\, \hat{B}_{V}\big)}\otimes\mathcal{O}_{V}\Big(f^{*}\big(4k-4)H\big)-\sum_{i=1}^{r}E_{i}\Big)\Big)%
$$
is surjective by the Shokurov vanishing theorem (see
Theorem~\ref{theorem:Shokurov}), where $\mathcal{L}(V,
\hat{B}_{V})$ is a subscheme of log canonical singularities of the
log pair $(V, \hat{B}_{V})$.

In the neighborhood of the point $\bar{P}$ the support of the
subscheme $\mathcal{L}(V, \hat{B}_{V})$ consists of the point
$\bar{P}$, which implies the existence of an effective divisor
$$
D\in \Big|f^{*}\Big(\big(4k-4\big)H\Big)-\sum_{i=1}^{r}E_{i}\Big|
$$
that does not pass through the point $\bar{P}$. Therefore, the
divisor $f(D)$ is a hypersurface of degree $4k-4$ that passes
through all points of the set $\mathrm{Sing}(V)\setminus P$ but
does not pass through the point $P$. We have $4k-4\leqslant 2n-5$,
which implies the existence of a hypersurface of degree $2n-5$
that contains the set $\mathrm{Sing}(V)\setminus P$ and does not
pass through the point $P$.

In general case we can apply the previous arguments to the linear
system $\hat{\mathcal{H}}$ instead of the linear system
$\mathcal{H}$, put $X=V$, and use the projective normality of
$V\subset\mathbb{P}^{4}$.
\end{proof}

\begin{corollary}
\label{corollary:factoriality-answer-to-Ciliberto}%
Suppose that the subset $\mathrm{Sing}(V)\subset\mathbb{P}^{4}$ is
a set-theoretical intersection hypersurfaces of degree $k<n/2$.
Then the hypersurface $V$ is $\mathbb{Q}$-factorial.
\end{corollary}

Every smooth surface on $V$ is a Cartier divisor when
$\mathrm{Sing}(V)<(n-1)^2$ due to \cite{CiGe03}, and it is natural
to expect that $V$ is $\mathbb{Q}$-factorial in the case when
$|\mathrm{Sing}(V)|<(n-1)^{2}$, which is proved only for
$n\leqslant 4$ (see \cite{FiWe89}, \cite{Ch04e}). The arguments of
the proof of
Proposition~\ref{proposition:factoriality-answer-to-Ciliberto} and
known properties of linear systems on rational surfaces are used
in \cite{Ch04t} to prove that the hypersurface $V$ is
$\mathbb{Q}$-factorial in the case when
$|\mathrm{Sing}(V)|\leqslant(n-1)^{2}/4$.

The main result of the given paper is the following theorem.

\begin{theorem}
\label{theorem:main}%
The following nodal threefolds are $\mathbb{Q}$-factorial:
\begin{itemize}
\item a complete intersection of the hypersurface $F$ and $G$ in
$\mathbb{P}^{5}$ of degree $n$ and $k$ respectively
such that $G$ smooth, $|\mathrm{Sing}(X)|\leqslant(n+k-2)(n-1)/5$, and  $n\geqslant k$;%
\item a double cover of a smooth hypersurface
$F\subset\mathbb{P}^{4}$ of degree $n\geqslant 2$ branched in a
surface $S\subset F$ that is cut out on $F$ by a hypersurface of
degree $2r\geqslant n$ such that the number of singular points of
the surface $S$ does not exceed $(2r+n-2)r/4$.
\end{itemize}
\end{theorem}

Nodal threefolds arise naturally in many problems of algebraic
geometry.

\begin{example}
\label{example:sextic-double-solid-bidegree-2-3}%
Let $Y$ be a general divisor of bi-degree $(2,3)$ in
$\mathbb{P}^1\times\mathbb{P}^3$ given by
$$
f_3\big(x,y,z,w\big)s^2+g_3\big(x,y,z,w\big)st+h_3\big(x,y,z,w\big)t^2=0,%
$$%
where $(s:t;x:y:z:w)$ are bihomogeneous coordinates, and  $f_3$,
$g_3$, $h_3$ are homogeneous polynomials of degree $3$. Let
$\xi:Y\to\mathbb{P}^3$ be a natural projection. Then $Y$ has $27$
rational curves $C_1, C_2, \cdots, C_{27}$ such that $-K_Y\cdot
C_i=0$, because the system of equations
$$
f_3\big(x,y,z,w\big)=g_3\big(x,y,z,w\big)=h_3\big(x,y,z,w\big)=0
$$
has exactly $27$ solutions. The projection $\xi$ has degree  $2$
outside of  $C_1, C_2, \cdots, C_{27}$, but
$$
X=\mathrm{Proj}\left(\bigoplus_{n\geqslant 0} H^0\Big(Y, \mathcal{O}_V\big(-nK_Y\big)\Big)\right)%
$$
is a double cover of $\mathbb{P}^3$ branched over a nodal surface
$$
g_3^2\big(x,y,z,w\big)-4f_3\big(x,y,z,w\big)h_3\big(x,y,z,w\big)=0,
$$
which implies that the threefold $X$ is nodal and has exactly $27$
ordinary double points that are images of the smooth rational
curves $C_1, C_2, \cdots, C_{27}$, which are contracted by the
morphism $\phi_{|-nK_Y|}:Y\to X$ for some natural $n\gg 0$. The
threefold $X$  is not $\mathbb{Q}$-factorial, and it is well
known, that the threefold $X$ is not rational (see \cite{Bar84},
\cite{Sob02}, \cite{Ch04d}).
\end{example}

The geometry of nodal threefolds is more complicated than of
smooth ones:
\begin{itemize}
\item every surface on smooth hypersurface in $\mathbb{P}^{4}$ is
a complete intersection due to Lefschetz theorem, which is no
longer true in the nodal case (see Example~\ref{example:Barth-sextic});%
\item the group of birational automorphisms of a smooth quartic
threefold is a finite group (see \cite{IsMa71}), which is no
longer true in the nodal case (see \cite{Pu88b}, \cite{Me03});%
\item smooth cubic threefolds are not rational (see
\cite{ClGr72}), but singular ones are rational.%
\end{itemize}

Isolated ordinary double point has two small resolutions, which
are birational via an ordinary flop (see \cite{We87},
\cite{Ko89}). Therefore, every nodal threefold having $k$ singular
points has exactly $2^{k}$ small resolutions, which all must be
non-projective in the $\mathbb{Q}$-factorial case, because every
exceptional curves must be homological to zero. Thus, it is quite
natural to expect that a singular nodal threefold is
$\mathbb{Q}$-factorial if and only if all its small resolutions
are not projective. The following example of L.\,Wotzlaw shows
that the latter is not true.

\begin{example}
\label{example:Lorenz} Let $\mathcal{I}_{5}$ be a quintic
hypersurface
$$
x_{5}-6x_{5}^{3}\sum_{i=0}^{6}x_{i}-27x_{5}\left(\Big(\sum_{i=0}^{5}x_{i}\Big)^{2}-4\sum_{i=0}^{5}\sum_{j=i+1}^{5}x_{i}x_{j}\right)-648x_{0}x_{1}x_{2}x_{3}x_{4}=0
$$
in
$\mathbb{P}^{5}\cong\mathrm{Proj}(\mathbb{C}[x_{0},x_{1},x_{2},x_{3},x_{4},x_{5}])$.
Then the quintic $\mathcal{I}_{5}$ is invariant under the standard
action of the Weil group $E_{6}$ on $\mathbb{P}^{5}$ by
reflection. Moreover, the quintic $\mathcal{I}_{5}$ is the only
invariant of degree $5$ of the Weil group $E_{6}$ under such
action (see \S 6 in \cite{Hu96}, \cite{Hu00}).

The singularities of the quintic $\mathcal{I}_{5}$ consist of
lines $L_{1},\ldots,L_{120}$, which intersect each other in points
$O_{1},\ldots,O_{36}$, the projectivization of a tangent cone to
$\mathcal{I}_{5}$ in $O_{k}$ is isomorphic to a so-called Segre
cubic (see \cite{Fi87}, \cite{Hu96}, \cite{Hu00}), but in every
point of the set
$$
\bigcup_{i=1}^{120} L_{i}\setminus \bigcup_{k=1}^{36} O_{k}
$$
the quintic $\mathcal{I}_{5}$ is locally isomorphic to a product
$\mathbb{C}\times\mathbb{A}_{1}$, where $\mathbb{A}_{1}$ is a
neighborhood of a three-dimensional ordinary double point.

Let $H_{\alpha}$ be a hyperplane section of the quintic
$\mathcal{I}_{5}$ that corresponds to a general point $\alpha$ of
the dual space $(\mathbb{P}^{5})^{*}$, and $T_{\beta}$ be a
hyperplane section of $\mathcal{I}_{5}$ that corresponds to a
general point
$\beta\in(\mathcal{I}_{5})^{*}\subset(\mathbb{P}^{5})^{*}$ and
tangents $\mathcal{I}_{5}$ in a point $P\in \mathcal{I}_{5}$.
Therefore, there is a five-dimensional family of hyperplane
sections $H_{\alpha}$, and four-dimensional family of tangent
hyperplane sections $H_{\beta}$. It follows from \cite{Hu96} that
both families are versal.

The variety $H_{\alpha}$ is a nodal hypersurface in
$\mathbb{P}^{4}$ of degree $5$ that has $120$ ordinary double
points $Q_{1},\ldots,Q_{120}$ such that $Q_{i}=L_{i}\cap
H_{\alpha}$. The variety $T_{\beta}$ is a nodal hypersurface of
degree $5$ that has $121$ ordinary double points $P_{1},\ldots,
P_{120}$ and $P$ such that $P_{i}=L_{i}\cap T_{\beta}$.

It follows from the Lefschetz theorem that
$\mathrm{rk}\,\mathrm{Pic}(H_{\alpha})=\mathrm{rk}\,\mathrm{Pic}(T_{\beta})=1$,
but it follows from \cite{Bo90} that
$\mathrm{rk}\,\mathrm{Cl}(H_{\alpha})=\mathrm{rk}\,\mathrm{Cl}(T_{\beta})=25$,
and $H_{\alpha}$ and $T_{\beta}$ are not $\mathbb{Q}$-factorial.

Let $\pi:\hat{T}_{\beta}\to T_{\beta}$ be a small resolution, and
$C_{i}$ and $C$ be the curves on $\hat{T}_{\beta}$ that are
contracted to the points $P_{i}$ and $P$ respectively. Then
$$
\mathcal{N}_{C\slash \hat{T}_{\beta}}\cong\mathcal{N}_{C_{i}\slash
\hat{T}_{\beta}}\cong\mathcal{O}_{\mathbb{P}^{1}}(-1)\oplus\mathcal{O}_{\mathbb{P}^{1}}(-1),
$$
where $C\cong C_{i}\cong\mathbb{P}^{1}$.

Let $\psi:\bar{H}_{\alpha}\to H_{\alpha}$ be a small resolution,
and $\tau:\hat{T}_{\beta}\to\bar{T}_{\beta}$ be a small
contraction of a smooth rational curve $C$ to an ordinary double
point $\bar{P}\in\bar{T}_{\beta}$. Then $\bar{P}$ is the only
singular point of the variety $\bar{T}_{\beta}$, and
five-dimensional family of smooth threefolds $\bar{H}_{\alpha}$ is
a smooth deformation of the threefold $\bar{T}_{\beta}$.
Therefore, there is an exact sequence (see \cite{We87})
$$
0\to H_3\big(\hat{T}_{\beta}, \mathbb{Z}\big)\to
H_3\big(\bar{T}_{\beta}, \mathbb{Z}\big) \to H_2\big(C,
\mathbb{Z}\big)\to H_2\big(\hat{T}_{\beta}, \mathbb{Z}\big)\to
H_2\big(\bar{T}_{\beta}, \mathbb{Z}\big)\to 0
$$
and an isomorphism $H_2(\bar{T}_{\beta}, \mathbb{Z})\cong
H_2(\bar{H}_{\alpha}, \mathbb{Z})$, but $$ h_2(\hat{T}_{\beta},
\mathbb{Z})=\mathrm{rk}\,\mathrm{Cl}(T_{\beta})=\mathrm{rk}\,\mathrm{Cl}(H_{\alpha})=h_2(\bar{H}_{\alpha},
\mathbb{Z}),
$$
which implies that the natural map $H_2(C, \mathbb{Z})\to
H_2(\hat{T}_{\beta}, \mathbb{Z})$ maps the whole homology group
$H_2(C, \mathbb{Z})$ to the zero. Hence, the curve $C$ is
homological to the zero on the smooth threefold $\hat{T}_{\beta}$,
which implies that $\hat{T}_{\beta}$ is not projective.
\end{example}

 Let us consider two examples, which are
inspired by the papers \cite{Kr00} and \cite{Me03}.

\begin{example}
\label{example:sextic-double-solid-Q-factoriality} Let $\pi:X\to
\mathbb{P}^3$ be the double cover ramified along a surface $S$
given by
$$
u^{2}+g^{2}_{3}\big(x,y,z,w\big)=h_{1}\big(x,y,z,w\big)f_{5}\big(x,y,z,w\big)\subset\mathbb{P}\big(1^{4},3\big)\cong\mathrm{Proj}\Big(\mathbb{C}[x,y,z,w,u]\Big),
$$
where $g_{3}$, $h_{1}$, and $f_{5}$ are general polynomials
defined over $\mathbb{R}$ of degree $3$, $1$, and $5$,
respectively. Then the double cover $X$ is not
$\mathbb{Q}$-factorial over $\mathbb{C}$ because the divisor
$h_{1}=0$ splits into two non-$\mathbb{Q}$-Cartier divisors
conjugated by $\mathrm{Gal}(\mathbb{C}\slash \mathbb{R})$ and
given by the equation
$$
\Big(u+\sqrt{-1}g_{3}\big(x,y,z,w\big)\Big)\Big(u-\sqrt{-1}g_{3}\big(x,y,z,w\big)\Big)=0.
$$

The sextic surface $S\subset\mathrm{Proj}(\mathbb{C}[x,y,z,w])$
has $15$ ordinary double points at the intersection points of the
three surfaces
$$
\big\{h_{1}(x,y,z,w)=0\big\}\cap\big\{g_{3}(x,y,z,w)=0\big\}\cap\big\{f_{5}(x,y,z,w)=0\big\},
$$
which gives $15$ simple double points of $X$. Introducing
variables  $s$ and $t$ defined by
$$
\left\{\aligned
&s={\frac{u+\sqrt{-1}g_{3}\big(x,y,z,w\big)} {h_{1}\big(x,y,z,w\big)}}={\frac{f_{5}\big(x,y,z,w\big)} {u-\sqrt{-1}g_{3}\big(x,y,z,w\big)}}\\
&t={\frac{u-\sqrt{-1}g_{3}\big(x,y,z,w\big)} {h_{1}\big(x,y,z,w\big)}}={\frac{f_{5}\big(x,y,z,w\big)} {u+\sqrt{-1}g_{3}\big(x,y,z,w\big)}}\\
\endaligned
\right.
$$
we can unproject $X$ in the sense of \cite{Re00} into two complete
intersections
$$
\left\{\aligned &V_{s}=\left\{\aligned
&s h_{1}\big(x,y,z,w\big)=u+\sqrt{-1}g_{3}\big(x,y,z,w\big)\\
&s \Big(u-\sqrt{-1}g_{3}\big(x,y,z,w\big)\Big)=f_{5}\big(x,y,z,w\big)\\
\endaligned
\right \}\subset \mathbb{P}\big(1^{4},3,2\big)\\
&V_{t}=\left\{\aligned
&t h_{1}\big(x,y,z,w\big)=u-\sqrt{-1}g_{3}\big(x,y,z,w\big)\\
&t \Big(u+\sqrt{-1}g_{3}\big(x,y,z,w\big)\Big)=f_{5}\big(x,y,z,w\big)\\
\endaligned
\right\}\subset \mathbb{P}\big(1^{4},3,2\big),\\
\endaligned
\right.
$$
respectively, which are not defined over $\mathbb{R}$. Eliminating
variable $u$, we get
$$
\left\{\aligned &V_{s}=\big\{s^2h_{1}-2\sqrt{-1}s g_{3}-f_{5}=0\big\}\subset \mathbb{P}(1^{4},2\big)\\
&V_{t}=\big\{t^2h_{1}+2\sqrt{-1}t g_{3}-f_{5}=0\big\}\subset \mathbb{P}\big(1^{4},2\big)\\
\endaligned
\right.
$$
and for the maps $\rho_{s}:X\dasharrow V_{s}$ and
$\rho_{t}:X\dasharrow V_{t}$ we obtain a commutative diagram
$$
\xymatrix{
&Y_s\ar[dl]_{\psi_s}  \ar[dr]^{\phi_s} && Y_t \ar[dl]_{\phi_t}\ar[dr]^{\psi_t}&\\
V_s && X\ar@{-->}[rr]^{\rho_t}\ar@{-->}[ll]_{\rho_s}&&V_t}
$$
with birational morphisms $\phi_{s}$, $\psi_{s}$, $\phi_{t}$, and
$\psi_{t}$ such that $\psi_{s}$ and $\psi_{t}$ are extremal
contractions in the sense of \cite{Co95}, while $\phi_{s}$ and
$\phi_{t}$ are flopping contractions.

The weighted hypersurfaces $V_{s}$ and $V_{t}$ are quasi-smooth
(see \cite{IF00}), which implies that they are
$\mathbb{Q}$-factorial and have Picard group $\mathbb{Z}$ (see
\cite{CalLy94}). The hypersurfaces $V_{s}$ and $V_{t}$ are
projectively isomorphic in $\mathbb{P}(1^{4},2)$ by the action of
$\mathrm{Gal}(\mathbb{C}\slash\mathbb{R})\cong\mathbb{Z}_{2}$. We
have
$$
\mathrm{Pic}\big(Y_{s}\big)\cong\mathrm{Pic}\big(Y_{t}\big)\cong \mathbb{Z}\oplus\mathbb{Z},%
$$
which gives $\mathrm{Cl}(X)=\mathbb{Z}\oplus\mathbb{Z}$. The
$\mathrm{Gal}(\mathbb{C}\slash\mathbb{R})$-invariant part of the
group $\mathrm{Cl}(X)$ is $\mathbb{Z}$, which implies that $X$ is
$\mathbb{Q}$-factorial over $\mathbb{R}$. The threefold $X$ is not
rational over $\mathbb{R}$ due to \cite{ChPa04}, but $X$ is also
not rational over $\mathbb{C}$ due to \cite{CPR}. Moreover, the
involution of $X$ interchanging fibers of $\pi$ induces a
non-biregular involution of $V_{s}$ which is regularized by
$\rho_{s}$ (see \cite{Ch04g}).
\end{example}

\begin{example}
\label{example:sextic-double-solid-quartic-with-double-point} Let
$V\subset \mathbb{P}^4$ be a general hypersurface of degree $4$
such that $V$ has exactly one ordinary double point $O$. Then $V$
is $\mathbb{Q}$-factorial and can be given by the equation
$$
t^{2}f_{2}\big(x,y,z,w\big)+tf_{3}\big(x,y,z,w\big)+f_{4}\big(x,y,z,w\big)=0\subset \mathbb{P}^4\cong\mathrm{Proj}\Big(\mathbb{C}[x,y,z,w,t]\Big),%
$$
where $O=(0:0:0:0:1)$. The threefold $V$ is known to be
non-ra\-ti\-o\-nal (see \cite{Me03}, \cite{Pu88b}), but the
projection $\phi:V\dasharrow \mathbb{P}^3$ from the singular point
$O$ has degree $2$ at a generic point of the threefold $V$ and
induces a non-biregular involution $\tau\in \mathrm{Bir}(V)$.

Let $f:Y\to V$ be the blow up of the point $O$. Then the linear
system $|-nK_{Y}|$ does not have base points for $n\gg 0$ and
gives a birational morphism $g:Y\to X$ contracting every curve
$C_i\subset Y$ such that $f(C_i)$ is a line on the quartic
threefold $V$ passing through the singular point $O$. We then
obtain the double cover $\pi:X\to\mathbb{P}^3$ ramified along the
nodal sextic surface $S\subset\mathbb{P}^3$ given by the equation
$$
f^2_{3}\big(x,y,z,w\big)-4f_{2}\big(x,y,z,w\big)f_{4}\big(x,y,z,w\big)=0.
$$

Each line $f(C_{i})$ corresponds to an intersection point of three
surfaces
$$
\big\{f_{2}(x,y,z,w)=0\big\}\cap\big\{f_{3}(x,y,z,w)=0\big\}\cap\big\{f_{4}(x,y,z,w)=0\big\}\subset\mathbb{P}^3\cong\mathrm{Proj}\Big(\mathbb{C}[x,y,z,w]\Big),
$$
which gives $24$ smooth rational curves $C_{1}, C_2, \cdots,
C_{24}$ such that
$$
\mathcal{N}_{Y\slash C_{i}}\cong \mathcal{O}_{C_{i}}(-1)\oplus
\mathcal{O}_{C_{i}}(-1)
$$
and $g$ is a  standard flopping contraction which maps every curve
$C_{i}$ to an ordinary double point of the threefold $X$. In
particular, the sextic $S\subset\mathbb{P}^3$ has exactly $24$
simple double points. However, the threefold $X$ is not
$\mathbb{Q}$-factorial and
$\mathrm{Cl}(X)=\mathbb{Z}\oplus\mathbb{Z}$.

Put $\rho=g\circ f^{-1}$. Then the involution
$\gamma=\rho\circ\tau\circ\rho^{-1}$ is biregular on $X$ and
interchanges the fibers of the double cover $\pi$. Thus the map
$\rho$ is a regularization of the birational non-biregular
involution $\tau$ in the sense of \cite{Ch04g}, while the
commutative diagram
$$
\xymatrix{
&Y\ar[dl]_{f}  \ar[dr]^{g} &&&& Y\ar[dl]_{g}\ar[dr]^{f}&\\
V\ar@{-->}[rr]^{\rho}&&X\ar@{->}[rr]_{\gamma}&&X&&\ar@{-->}[ll]_{\rho}V}
$$%
is a decomposition of $\tau$ in a sequence of so-called Sarkisov
links (see \cite{Co95}, \cite{CPR}, \cite{Is96b}).

Suppose that $f_{2}(x,y,z,w)$ and $f_{4}(x,y,z,w)$ are defined
over $\mathbb{Q}$ and
$$
f_{3}\big(x,y,z,w\big)=\sqrt{2}g_{3}\big(x,y,z,w\big),
$$
where $g_{3}(x,y,z,t)$ is defined over $\mathbb{Q}$. Then the
threefold $V$ is defined over $\mathbb{Q}(\sqrt{2})$, but the
hypersurface $V$ is not defined over $\mathbb{Q}$, because the
threefold $V$ is not invariant under the action of
$\mathrm{Gal}(\mathbb{Q}(\sqrt{2})\slash\mathbb{Q})$. However, the
sextic surface $S\subset\mathbb{P}^3$ is given by the equation
$$
2g^2_{3}\big(x,y,z,w\big)-4f_{2}\big(x,y,z,w\big)f_{4}\big(x,y,z,w\big)=0\subset
\mathbb{P}^3\cong\mathrm{Proj}\Big(\mathbb{Q}[x,y,z,w]\Big),
$$
which implies that $X$ is defined over $\mathbb{Q}$ as well.
Moreover, the
$\mathrm{Gal}(\mathbb{Q}(\sqrt{2})\slash\mathbb{Q})$-invariant
part of the group $\mathrm{Cl}(X)$ is $\mathbb{Z}$, which implies
that $X$ is $\mathbb{Q}$-factorial over $\mathbb{Q}$.
\end{example}

Thus, the condition of $\mathbb{Q}$-factoriality depends also on
the field of definition.

\vskip 0.5cm

The author would like to cordially thank M.\,Gri\-nen\-ko,
V.\,Iskov\-s\-kikh, J.\,Park, Yu.\,Pro\-kho\-rov, V.\,Sho\-ku\-rov
and L.\,Wotzlaw for helpful conversations.

\section{Preliminaries.}
\label{section:preliminaries}

The following result is well known (see \cite{Cl83}, \cite{Sch85},
\cite{We87}, \cite{Di90}, \cite{Cy01}).

\begin{theorem}
\label{theorem:Cynk} Let $W$ be a smooth fourfold, $Y$ be an ample
reduced and irreducible divisor on the fourfold $W$ such that the
only singularities of the threefold $Y$ are nodal and
$$
h^{2}\big(\Omega^{1}_{W}\big)=h^{3}\Big(\Omega^{1}_{W}\otimes \mathcal{O}_{W}\big(-Y\big)\Big)=h^{1}\big(\mathcal{O}_{W}\big)=h^{2}\big(\mathcal{O}_{W}\big)=0,%
$$
and let ${\tilde Y}$ be a small resolution of the threefold $Y$.
Then
\begin{multline*}
h^{1}\big(\mathcal{O}_{\tilde{Y}}\big)=h^{2}\big(\mathcal{O}_{\tilde{Y}}\big)=0,\ h^{1}\big(\Omega^{1}_{\tilde Y}\big)=h^{1}\big(\Omega^{1}_{W}\big)+\delta,\\
h^{2}\big(\Omega^{1}_{\tilde Y}\big)=h^{0}\Big(K_{W}\otimes \mathcal{O}_{W}\big(2Y\big)\Big)+h^{3}\big(\mathcal{O}_{W}\big)-h^{0}\Big(K_{W}\otimes \mathcal{O}_{W}\big(Y\big)\Big)-\\ %
-h^{3}\big(\Omega^{1}_{W}\big)-h^{4}\Big(\Omega^{1}_{W}\otimes \mathcal{O}_{W}\big(-Y\big)\Big)-\big|\mathrm{Sing}(Y)\big|+\delta,%
\end{multline*}
where $\delta$ is the number of dependent conditions that
vanishing at the nodes of the threefold $Y$ imposes on the global
sections of the line bundle $K_{W}\otimes \mathcal{O}_{W}(2Y)$.
\end{theorem}

The proof of Theorem~\ref{theorem:Cynk} in \cite{Cy01} implies the
following result.

\begin{corollary}
\label{corollary:Cynk} Let $W$ be a smooth fourfold, and $Y$ be an
ample reduced and irreducible divisor on the fourfold $W$ such
that the threefold $Y$ is nodal. Suppose that
$$
h^{2}\big(\Omega^{1}_{W}\big)=h^{1}\big(\mathcal{O}_{W}\big)=h^{2}\big(\mathcal{O}_{W}\big)=0,%
$$
but singular points of the threefold $Y$ impose independent linear
conditions on the global sections of the line bundle
$K_{W}\otimes\mathcal{O}_{W}(2Y)$. Then $Y$ is
$\mathbb{Q}$-factorial.
\end{corollary}

The following result is due to \cite{Bes83}.

\begin{theorem}
\label{theorem:Bese} Let $\pi:Y\to\mathbb{P}^2$ be a blow up of
points $P_{1},\ldots, P_{s}$ such that
$$
s\leqslant{\frac{d^{2}+9d+10}{6}}
$$
and at most $k(d+3-k)-2$ points among the points $P_{1},\ldots,
P_{s}$ are contained in a curve of degree $k\leqslant (d+3)/2$ for
some natural number $d\geqslant 3$. Then the linear system
$$
\Big|\pi^{*}\Big(\mathcal{O}_{\mathbb{P}^2}\big(d\big)\Big)-\sum_{i=1}^{s}E_{i}\Big|
$$
does not have base points, where $E_{i}$ is the $\pi$-exceptional
divisor such that $\pi(E_{i})=P_{i}$.
\end{theorem}

In the case $d=3$ the claim of Theorem~\ref{theorem:Bese} is a
base point freeness of the anticanonical linear system of a
\emph{weak del Pezzo surface} of degree $9-s\geqslant 2$ (see
\cite{De80}, \cite{HiWa81}, \cite{Mae94}).

\begin{corollary}
\label{corollary:Bese} Let $\Sigma$ be a finite subset of
$\mathbb{P}^{2}$ and $d\geqslant 3$ be a natural number such that
$$
\big|\Sigma\big|\leqslant {\frac{d^{2}+9d+16}{6}}
$$
and at most $k(d+3-k)-2$ points of the set $\Sigma$ lie on a
possibly reducible plane curve of degree $k\leqslant(d+3)/2$. Then
for every point $P\in\Sigma$ there is a curve on $\mathbb{P}^{2}$
of degree $d$ that passes through all points of the set
$\Sigma\setminus P$ and does not pass through the point $P$.
\end{corollary}

The claim of Theorem~\ref{theorem:Bese} is strengthen in
\cite{DaGe88} in the following way.

\begin{theorem}
\label{theorem:Davis-Geramita} Let $\pi:Y\to\mathbb{P}^2$ be a
blow up of points $P_{1},\ldots, P_{s}$ such that
$$
s\leqslant \mathrm{max}\Big\{\Big\lfloor(d+3)/2\Big\rfloor\Big(d+3-\Big\lfloor(d+3)/2\Big\rfloor\Big)-1, \Big\lfloor(d+3)/2\Big\rfloor^{2}\Big\},%
$$
and at most $k(d+3-k)-2$ points among the points $P_{1},\ldots,
P_{s}$ are contained in a curve of degree $k\leqslant (d+3)/2$ for
some natural number $d\geqslant 3$. Then the linear system
$$
\Big|\pi^{*}\Big(\mathcal{O}_{\mathbb{P}^2}\big(d\big)\Big)-\sum_{i=1}^{s}E_{i}\Big|
$$
does not have base points, where $E_{i}$ is the $\pi$-exceptional
divisor such that $\pi(E_{i})=P_{i}$.
\end{theorem}

\section{Connectedness principle.}
\label{section:connectedness}

Let $X$ be a smooth variety, and $B_{X}=\sum_{i=1}^{k}a_{i}B_{i}$
be a $\mathbb{Q}$-divisor, where $B_{i}$ is a prime divisor and
$a_{i}$ is a positive rational number. Let $\pi:Y\to X$ be a
bi\-ra\-ti\-onal morphism such that $Y$ is smooth, and the union
of all the  proper transforms of the divisors $B_{i}$ and all the
$\pi$-exceptional divisors form a divisor with simple normal
crossing. Let $B_{Y}$ be the proper transform of $B_{X}$ on the
variety $Y$, and put
$$
B^{Y}=B_{Y}-\sum_{i=1}^{n}c_{i}E_{i},
$$
where $E_{i}$ is an $\pi$-ex\-cep\-ti\-o\-nal divisor and $c_{i}$
is a rational number such that the equivalence
$$
K_{Y}+B^{Y}\sim_{\mathbb{Q}}\pi^{*}\big(K_{X}+B_{X}\big)
$$
holds. Then the log pair $(Y, B^{Y})$ is called the log pull back
of the log pair $(X, B_{X})$ with respect to the birational
morphism $\pi$, while the number $c_{i}$ is called the discrepancy
of the log pair $(X, B_{X})$ in the $\pi$-exceptional divisor
$E_i$.

\begin{definition}
\label{definition:log-center} A proper irreducible subvariety
$Z\subset X$ is called a center of log canonical singularities of
the log pair $(X, B_{X})$ if there is a divisor $E$ on $Y$
contained in the support of the effective part of the divisor
$\lfloor B^{Y}\rfloor$ such that $\pi(E)=Z$.
\end{definition}

In particular, the proper irreducible subvariety
$\pi(E_{i})\subset X$ is a center of log canonical singularities
of the log pair $(X, B_{X})$ if $c_{i}\leqslant -1$. Similarly,
the prime divisor $B_{i}$ is center of log canonical singularities
of the log pair $(X, B_{X})$ if $a_{i}\geqslant 1$.

The set of all centers of log canonical singularities of the log
pair $(X, B_{X})$ are usually denoted as $\mathbb{LCS}(X, B_{X})$.
Similarly, the union of all centers of log canonical singularities
of the log pair $(X, B_{X})$ considered as a proper subset of the
variety $X$ are called the locus of log canonical singularities of
the log pair $(X, B_{X})$ and denoted as $LCS(X, B_{X})$.

\begin{example}
\label{example:smooth-point} Let $O$ be a smooth point on $X$.
Then the inequality
$\mathrm{mult}_{O}(B_{X})\geqslant\mathrm{dim}(X)$ implies that
$O\in\mathbb{LCS}(X, B_{X})$. Moreover, the inequality
$\mathrm{mult}_{O}(B_{X})\geqslant 1$ holds in the case when
$O\in\mathbb{LCS}(X, B_{X})$ and the boundary $B_{X}$ is
effective.
\end{example}

\begin{remark}
\label{remark:log-canonical-reduction} Let $H$ be a general
hyperplane section of the variety $X$, and $Z$ be subvariety of
the variety $X$ that is an element of the set $\mathbb{LCS}(X,
B_{X})$. Then every component of the intersection $Z\cap H$ is
contained in $\mathbb{LCS}(H, B_{X}\vert_{H})$.
\end{remark}

\begin{example}
\label{example:smooth-point-and-log-pull-back} Let $O$ be a smooth
point of the variety $X$. Suppose that $O$ is a center of log
canonical singularities of the log pair $(X, B_{X})$. Let $f:V\to
X$ be the blow up of the point $O$, and $E$ be the $f$-exceptional
divisor. Then either $E\in\mathbb{LCS}(V, B^{V})$, or there is a
proper irreducible subvariety $Z\subset E$ that is a center of log
canonical singularities of the log pair $(V, B^{V})$. Moreover,
the exceptional divisor $E$ is a center of log canonical
singularities of the log pair $(V, B^{V})$ if and only if
$\mathrm{mult}_{O}(B_{X})\geqslant\mathrm{dim}(X)$.
\end{example}

\begin{definition}
\label{definition:log-subscheme} The subscheme associated to the
ideal sheaf
$$
\mathcal{I}\big(X, B_{X}\big)=f_{*}\Big(\mathcal{O}_{Y}\big(\lceil -B^{Y}\rceil\big)\Big)%
$$
is called the log canonical singularity subscheme of $(X, B_{X})$
and denoted as $\mathcal{L}(X, B_{X})$.
\end{definition}

The support of the subscheme $\mathcal{L}(X, B_{X})$ consists of
the set-theoretic union of all centers of log canonical
singularities of the log pair $(X, B_{X})$, which implies that
$$
\mathrm{Supp}\Big(\mathcal{L}\big(X, B_{X}\big)\Big)=LCS\big(X, B_{X}\big)\subset X.%
$$

The following result is the Shokurov vanishing theorem (see
\cite{Sho93}, \cite{Ko91}, \cite{Ko97}, \cite{Am99}).

\begin{theorem}
\label{theorem:Shokurov} Suppose that $K_{X}+B_{X}+H$ is
numerically equivalent to a Cartier divisor, where $H$ a
$\mathbb{Q}$-divisor on the variety $X$ that is nef and
big\footnote{It should be pointed out that a $\mathbb{Q}$-Cartier
divisor $H\in\mathrm{Div}(X)\otimes\mathbb{Q}$ is called
numerically effective or nef if for every curve $C\subset X$ the
inequality $H\cdot C\geqslant 0$ holds. A numerically effective
divisor $H$ is called big if the inequality $H^{n}>0$ holds, where
$n=\mathrm{dim}(X)$.}. Then for every $i>0$ we have
$$
H^{i}\Big(X, \mathcal{I}\big(X, B_{X}\big)\otimes
\mathcal{O}_{X}\big(K_{X}+B_{X}+H\big)\Big)=0.
$$
\end{theorem}

\begin{proof}
It follows from the Kawamata-Viehweg vanishing theorem (see
\cite{Ka82}, \cite{Vi82}) that
$$
R^{i}f_{*}\Big(f^{*}\big(K_{X}+B_{X}+H\big)+\lceil -B^{Y}\rceil\Big)=0%
$$
for all $i>0$. The degeneration of the local--to--global spectral
sequence and
$$
R^{0}f_{*}\Big(f^{*}\big(K_{X}+B_{X}+H\big)+\lceil -B^{Y}\rceil\Big)=\mathcal{I}\big(X, B_{X}\big)\otimes \mathcal{O}_{X}\big(K_{X}+B_{X}+H\big)\Big)%
$$
imply that
$$
H^{i}\Big(X, \mathcal{I}\big(X, B_{X}\big)\otimes\mathcal{O}_{X}\big(K_{X}+B_{X}+H\big)\Big)=H^{i}\Big(Y, f^{*}\big(K_{X}+B_{X}+H\big)+\lceil -B^{Y}\rceil\Big)%
$$
for $i\geqslant 0$. On the other hand, we have
$$
H^{i}\Big(Y, f^{*}\big(K_{X}+B_{X}+H\big)+\lceil -B^{Y}\rceil\Big)=0%
$$
for $i>0$ by the Kawamata-Viehweg vanishing theorem.
\end{proof}

The claim of Theorem~\ref{theorem:Shokurov} implies the following
result.

\begin{lemma}\label{lemma:non-vanishing}
Let $\mathcal{M}$ be a linear subsystem in
$|\mathcal{O}_{\mathbb{P}^n}(k)|$ such that the base locus of the
linear system $\mathcal{M}$ is zero-dimensional. Then the points
of the base locus of $\mathcal{M}$ impose independent linear
conditions on homogeneous forms on $\mathbb{P}^n$ of degree
$n(k-1)$.
\end{lemma}

\begin{proof}
Let $\Lambda$ be the base locus of the linear system
$\mathcal{M}$, and $H_1, \cdots, H_r$ be general divisors in the
linear system $\mathcal{M}$, where $r$ is sufficiently big. Put
$B_{\mathbb{P}^{n}}=\frac{n}{r}\sum_{i=1}^r H_i$. Then the
singularities of the log pair $(\mathbb{P}^n, B_{\mathbb{P}^{n}})$
are log terminal (see \cite{Ko97}) outside of the set $\Lambda$,
but
$$
\mathrm{mult}_{P}\Big(B_{\mathbb{P}^{n}}\Big)=n\sum_{i=1}^r \frac{\mathrm{mult}_{P}\big(H_i\big)}{r}\geqslant n%
$$
for every point $P\in\Lambda$. Thus, we have
$\mathrm{Supp}(\mathcal{L}(\mathbb{P}^n,
B_{\mathbb{P}^{n}}))=\Lambda$.

Since
$K_{\mathbb{P}^n}+B_{\mathbb{P}^{n}}+H\sim_{\mathbb{Q}}n(k-1)H$,
where $H$ is a hyperplane in $\mathbb{P}^{n}$, we see that
$$
H^1\left(\mathbb{P}^n, \mathcal{I}\Big(\mathbb{P}^n,
B_{\mathbb{P}^{n}}\Big)\otimes\mathcal{O}_{\mathbb{P}^n}\Big(n\big(k-1\big)\Big)\right)=0
$$
by Theorem~\ref{theorem:Shokurov}. Hence, the points of $\Lambda$
impose independent linear conditions on homogeneous forms of
degree $n(k-1)$, because $\mathrm{Supp}(\mathcal{L}(\mathbb{P}^n,
B_{\mathbb{P}^{n}}))=\Lambda$.
\end{proof}

\section{Complete intersections.}
\label{section:complete-intersections}

Let $X$ be a complete intersection of hypersurfaces $F$ and $G$ in
$\mathbb{P}^{5}$ such that the singularities of $X$ are nodal. Put
$n=\mathrm{deg}(F)$ and $k=\mathrm{deg}(G)$. Suppose that
$n\geqslant k$.

\begin{example}
\label{example:factoriality-nodal-complete-intersections-non-Q-factoriality}
Let $F$ and $G$ be general hypersurfaces containing a plane. Then
$X$ is nodal and not $\mathbb{Q}$-factorial, both $F$ and $G$ are
smooth, and $|\mathrm{Sing}(X)|=(n+k-2)^{2}$.
\end{example}

The following result is proved in \cite{CiGe04}.

\begin{theorem}
\label{theorem:factoriality-nodal-complete-intersections-Ciliberto}%
Suppose that $G$ is smooth and $|\mathrm{Sing}(X)|\leqslant 3n/8$.
Then $X$ is $\mathbb{Q}$-factorial.
\end{theorem}

In this section we prove the following result.

\begin{theorem}
\label{theorem:factoriality-nodal-complete-intersections-main}%
Suppose that $G$ is smooth. Then $X$ is $\mathbb{Q}$-factorial in
the case when
$$
\big|\mathrm{Sing}(X)\big|\leqslant\frac{\big(n+k-2\big)\big(n-1\big)}{5}.
$$
\end{theorem}

The claim of
Theorem~\ref{theorem:factoriality-nodal-complete-intersections-main}
is not true in the case when the hypersurface $G$ is singular.

\begin{example}
\label{example:non-Q-factoriality-singular} Let
$Q\subset\mathbb{P}^{5}$ be a smooth quadric surface, and $G$ be a
cone over the quadric surface $Q$ whose vertex is a general line
$L\subset\mathbb{P}^{5}$. Take a general hypersurface
$F\subset\mathbb{P}^{5}$ of degree $n$. Let $X$ be the complete
intersection of the hypersurfaces $G$ and $F$. Then $X$ is a nodal
threefold of degree $2n$ and $|\mathrm{Sing}(X)|=n$. Let $\Omega$
be a linear subspace in $\mathbb{P}^{5}$ spanned by a line
contained in $Q$ and a line $L$. Then $\Omega\subset G$, the
surface $\Omega\cap F$ has degree $n$ and is not a
$\mathbb{Q}$-Cartier divisor on the threefold $X$.
\end{example}

In the case $k=1$ the claim of
Theorem~\ref{theorem:factoriality-nodal-complete-intersections-main}
follows from \cite{Ch04t}.

\begin{conjecture}
\label{conjecture:factoriality-nodal-complete-intersections-factoriality}
Suppose that $G$ is smooth. Then $X$ is $\mathbb{Q}$-factorial in
the case when
$$
\big|\mathrm{Sing}(X)\big|\leqslant(n+k-2)^{2}.
$$
\end{conjecture}

Suppose that $G$ is smooth. Then the following result follows from
Corollary~\ref{corollary:Cynk}.

\begin{proposition}
\label{proposition:factoriality-nodal-complete-intersections-defect}%
The threefold $X$ is $\mathbb{Q}$-factorial in the case when its
singular points impose independent linear conditions on the
sections in
$H^{0}(\mathcal{O}_{\mathbb{P}^{5}}(2n+k-6)\vert_{G})$.
\end{proposition}

\begin{corollary}
\label{corollary:factoriality-nodal-complete-intersections-defect}%
Suppose that $|\mathrm{Sing}(X)|\leqslant 2n+k-5$. Then $X$ is
$\mathbb{Q}$-factorial.
\end{corollary}

The variety $X$ is $\mathbb{Q}$-factorial if and only if the group
$\mathrm{Cl}(X)$ is generated by the class of a hyperplane section
(see Remark~\ref{remark:factoriality-of-hypersurface}). Every
surface contained in the threefold $X$ is a complete intersection
in $\mathbb{P}^{5}$ in the case when $X$ is
$\mathbb{Q}$-factorial.

Now we prove
Theorem~\ref{theorem:factoriality-nodal-complete-intersections-main}.
Suppose thata $|\mathrm{Sing}(X)|\leqslant(n+k-2)(n-1)/5$, but the
hypersurface $G$ is smooth. We have $n=\mathrm{deg}(F)\geqslant
k=\mathrm{deg}(G)$. Let us show that the singular points of the
complete intersection $X\subset\mathbb{P}^{5}$ impose independent
linear conditions on the hypersurface in $\mathbb{P}^{5}$ of
degree $2n+k-6$, which implies the claim of
Theorem~\ref{theorem:factoriality-nodal-complete-intersections-main}.

The claim of
Theorem~\ref{theorem:factoriality-nodal-complete-intersections-main}
follows from \cite{Ch04t} in the case $k=1$, and in the case
$4\geqslant n$ the claim of
Theorem~\ref{theorem:factoriality-nodal-complete-intersections-main}
follows
Corollary~\ref{corollary:factoriality-nodal-complete-intersections-defect}.
Thus, we assume that $k\geqslant 2$ and $n\geqslant 5$

\begin{lemma}
\label{lemma:factoriality-nodal-complete-intersections-isolated-singularities}
There is a hypersurface $\hat{F}\subset\mathbb{P}^{5}$ of degree
$n$ such that the threefold $X$ is a complete intersection of the
hypersurfaces $\hat{F}$ and $G$, but
$\mathrm{Sing}(\hat{F})\subseteq \mathrm{Sing}(X)$.
\end{lemma}

\begin{proof}
The threefold $X$ is given by the system of equations
$$
\left\{\aligned
&f(x_0,x_1,x_2,x_3,x_{4},x_{5})=0\\
&g(x_0,x_1,x_2,x_3,x_{4},x_{5})=0\\
\endaligned
\right.
\subset\mathbb{P}^{5}\cong\mathrm{Proj}\Big(\mathbb{C}[x_0,x_1,x_2,x_3,x_{4},x_{5}]\Big),
$$
where $f$ and $g$ be are homogeneous polynomials of degree $n$ and
$k$  that define the hypersurface $F$ and $G$ respectively.
Consider linear system
$$
\mathcal{L}=\big|\lambda
f+h(x_0,x_1,x_2,x_3,x_{4},x_{5})g\big|\subset\big|\mathcal{O}_{\mathbb{P}^5}(n)\big|,
$$
where $\lambda\in\mathbb{C}$, and $h$ is a homogeneous polynomial
of degree $n-k$. Then the base locus of the linear system
$\mathcal{L}$ is the variety $X$. The Bertini theorem implies the
existence of a hypersurface $\hat{F}\subset\mathcal{L}$ such that
$X=\hat{F}\cap G$, but $\mathrm{Sing}(\hat{F})\subseteq
\mathrm{Sing}(X)$.
\end{proof}

We may assume that $\mathrm{Sing}(F)\subseteq\mathrm{Sing}(X)$.

\begin{definition}
\label{definition:general-position-dva} We say that the points of
a subset $\Gamma\subset{\mathbb P}^{r}$ have property $\bigstar$
in the case when at most $t(n+k-2)$ points of the set $\Gamma$ lie
on a curve in $\mathbb{P}^{r}$ of degree $t\in\mathbb{N}$.
\end{definition}

Let $\Sigma=\mathrm{Sing}(X)\subset\mathbb{P}^{5}$.

\begin{proposition}
\label{proposition:factoriality-nodal-complete-intersections-nodes-in-general-position}
The points of the subset $\Sigma\subset\mathbb{P}^{5}$ have
property $\bigstar$.
\end{proposition}

\begin{proof}
The hypersurface $F\subset\mathbb{P}^{5}$ can be given by the
equation
$$
f\big(x_0,x_1,x_2,x_3,x_4,x_5\big)=0\subset\mathbb{P}^{5}\cong\mathrm{Proj}\Big(\mathbb{C}[x_0,x_1,x_2,x_3,x_4,x_5]\Big),
$$
where $f$ is a homogeneous polynomial of degree $n$, and $G$ can
be given by the equation
$$
g\big(x_0,x_1,x_2,x_3,x_4,x_5\big)=0\subset\mathbb{P}^{5}\cong\mathrm{Proj}\Big(\mathbb{C}[x_0,x_1,x_2,x_3,x_4,x_5]\Big),
$$
where $g$ is a homogeneous polynomial of degree $k$. Then the set
$\Sigma$ is given by the vanishing of polynomials $f$ and $g$, and
by vanishing of all minors of size $1$ of the matrix
$$
\left(
\begin{matrix}
\frac{\partial f}{\partial x_0} & \frac{\partial f}{\partial x_1} &\frac{\partial f}{\partial x_2} & \frac{\partial f}{\partial x_3} & \frac{\partial f}{\partial x_4} & \frac{\partial f}{\partial x_5}\cr%
\frac{\partial g}{\partial x_0} & \frac{\partial g}{\partial x_1} &\frac{\partial g}{\partial x_2} & \frac{\partial g}{\partial x_3} & \frac{\partial g}{\partial x_4} & \frac{\partial g}{\partial x_5}\cr%
\end{matrix}\right),
$$
which implies that $\Sigma$ is a set-theoretical intersection of
hypersurfaces of degree $n+k-2$, which concludes the proof.
\end{proof}

Take an arbitrary point $P\in\Sigma$. Then we must show that there
is a hypersurface of degree $2n+k-6$ that contains the set
$\Sigma\setminus P$ and does not contain the point $P$.

\begin{lemma}
\label{lemma:factoriality-nodal-complete-intersections-complete-intersection-I}
Suppose that there is a plane $\Pi\subset\mathbb{P}^{5}$ such that
$\Sigma\subset\Pi\subset\mathbb{P}^{5}$. Then there is a
hypersurface of degree $2n+k-6$ that contains $\Sigma\setminus P$
and does not contain $P$.
\end{lemma}

\begin{proof}
We want to apply Corollary~\ref{corollary:Bese} to
$\Sigma\subset\Pi$ and $d=2n+k-6\geqslant 6$. Let us check that
all conditions of Corollary~\ref{corollary:Bese} are satisfied.

Suppose that $|\Sigma|>(d^{2}+9d+16)/6$. Then
$$
{\frac{(n+k-2)(n-1)}{5}}>{\frac{(2n+k-6)^{2}+9(2n+k-6)+16}{6}},
$$
where $n\geqslant 5$ and $k\geqslant 2$. Put $A=n+k\geqslant 7$.
Then
$$
0>\big(A+n-6\big)^{2}+9\big(A+n-6\big)+16-6An=5A^{2}-3A-10+5n^{2}-3n+4An\geqslant
464,
$$
which is a contradiction.

Now must show that at most $t(2n+k-3-t)-2$ points of the set
$\Sigma$ lie on a curve of degree $t\leqslant (2n+k-3)/2$.
However, at most $t(n+k-2)$ points of the set $\Sigma$ lie on a
curve of degree $t$ by
Proposition~\ref{proposition:factoriality-nodal-complete-intersections-nodes-in-general-position}.
In particular, in the case  $t=1$ we have
$$
t\big(2n+k-3-t\big)-2=2n+k-6\geqslant n+k-2=t\big(n-1\big),
$$
because $n\geqslant 5$. In the case when $t>1$ it is enough to
show that
$$
t\big(2n+k-3-t\big)-2\geqslant t\big(n+k-2\big)
$$
for every $t\leqslant (2n+k-3)/2$ such that
$t(2n+k-3-t)-2<|\Sigma|$. We have
$$
t(2n+k-3-t)-2\geqslant t(n+k-2)\iff n-1>t
$$
in the case when $t>1$. Therefore, we may assume that $t\geqslant
n-1$, which implies that
$$
t\big(2n+k-3-t\big)-2\geqslant
\big(n-1\big)\big(n+k-2\big)>\big|\Sigma\big|.
$$

Therefore, it follows from Corollary~\ref{corollary:Bese} that
there is a curve $C\subset\Pi$ of degree $2n+k-6$ that contains
the set $\Sigma\setminus P$ and does not contains $P$. Let $Y$ be
a general four-dimensional cone in $\mathbb{P}^{5}$ over the curve
$C$. Then $Y$ is the required hypersurface.
\end{proof}

Let $\Pi$ and $\Gamma$ be sufficiently general planes in
$\mathbb{P}^{5}$, and $\psi:\mathbb{P}^{5}\dasharrow\Pi$ be a
projection from the plane $\Gamma$. Put
$\Sigma^{\prime}=\psi(\Sigma)\subset\Pi\cong\mathbb{P}^{2}$ and
$\hat{P}=\psi(P)\in\Sigma^{\prime}$.

\begin{lemma}
\label{lemma:factoriality-nodal-complete-intersections-complete-intersection-II}
Suppose that the points of $\Sigma^{\prime}\subset\Pi$ have
property $\bigstar$. Then there is a hy\-per\-sur\-face of degree
$2n+k-6$ that contains $\Sigma\setminus P$ and does not contain
$P$.
\end{lemma}

\begin{proof}
The proof of
Lemma~\ref{lemma:factoriality-nodal-complete-intersections-complete-intersection-I}
implies the existence of a curve $C\subset\Pi$ of degree $2n+k-6$
that contains $\Sigma^{\prime}\setminus\hat{P}$ but does not pass
through the point $\hat{P}$. Let $Y\subset\mathbb{P}^{5}$ be the
cone over the curve $C$ whose vertex is $\Gamma$. Then $Y$ is a
hypersurface in $\mathbb{P}^{5}$ of degree $2n+k-6$ that passes
through all points of the set $\Sigma\setminus P$ and does not
contain the point $P\in\Sigma$.
\end{proof}

Therefore, we may assume that the points of the set
$\Sigma^{\prime}\subset\Pi\cong\mathbb{P}^{2}$ does not have the
property $\bigstar$. There is subset
$\Lambda_{r}^{1}\subset\Sigma$ such that
$|\Lambda_{r}^{1}|>r(n+k-2)$, but the subset
$$
\tilde{\Lambda}_{r}^{1}=\psi(\Lambda_{r}^{1})\subset\Sigma^{\prime}\subset\Pi\cong\mathbb{P}^{2}
$$
is contained in a curve $C\subset\Pi$ of degree $r$. Moreover, we
may assume that $r$ is the smallest natural number having such
property, which implies that the curve $C$ is irreducible and
reduced. We can iterate the construction of $\Lambda_{r}^{1}$ to
get the disjoint union of subsets
$$
\bigcup_{j=r}^{l}\bigcup_{i=1}^{c_{j}}\Lambda_{j}^{i}\subset\Sigma
$$
such that $|\Lambda_{j}^{i}|>j(n+k-2)$, the points of the set
$$
\tilde{\Lambda}_{j}^{i}=\psi(\Lambda_{j}^{i})\subset\Sigma^{\prime}
$$
lie on an irreducible curve in $\Pi\cong\mathbb{P}^{2}$ of degree
$j$, and the points of the subset
$$
\bar{\Sigma}=\Sigma^{\prime}\setminus\bigcup_{j=r}^{l}\bigcup_{i=1}^{c_{j}}\tilde{\Lambda}_{j}^{i}\subsetneq\Sigma^{\prime}\subset\Pi\cong\mathbb{P}^{2}
$$
have property $\bigstar$, where $c_{j}\geqslant 0$. Then
$c_{r}>0$ and
\begin{equation}
\label{equation:factoriality-nodal-complete-intersections-number-of-good-points}
\big|\bar{\Sigma}\big|<{\frac{\big(n+k-2\big)\big(n-1\big)}{5}}-\sum_{i=r}^{l}c_{i}\big(n-1\big)i={\frac{n+k-2}{5}}\Big(n-1-\sum_{i=r}^{l}5ic_{i}\Big).
\end{equation}

\begin{corollary}
\label{corollary:factoriality-nodal-complete-intersections-from-the-number-of-good-points}
The inequality $\sum_{i=r}^{l}ic_{i}<(n-1)/5$ holds.
\end{corollary}

In particular, we have $j<(n-1)/5$ in the case when
$\Lambda_{j}^{i}\ne\varnothing$.

\begin{lemma}
\label{lemma:factoriality-nodal-complete-intersections-complete-intersection-III}
Suppose that $\Lambda_{j}^{i}\ne\varnothing$. Let $\mathcal{M}$ be
a linear system of hypersurfaces in $\mathbb{P}^{5}$ of degree $j$
that contains $\Lambda_{j}^{i}$. Then the base locus of the linear
system $\mathcal{M}$ is zero-dimensional.
\end{lemma}

\begin{proof}
The construction of the set $\Lambda_{j}^{i}$ implies that all
points of the subset
$$
\tilde{\Lambda}_{j}^{i}=\psi\big(\Lambda_{j}^{i}\big)\subset\Sigma^{\prime}\subset\Pi\cong\mathbb{P}^{2}
$$
are contained in an irreducible curve $C\subset\Pi$ of degree $j$.
Let $Y$ be a cone in $\mathbb{P}^{5}$ over the curve $C$ whose
vertex is some plane $\Upsilon\subset\mathbb{P}^{5}$. Then $Y$ is
a hypersurface in $\mathbb{P}^{5}$ of degree $j$ that contains all
points of the set $\Lambda_{j}^{i}$, which implies that
$Y\in\mathcal{M}$.

Suppose that the base locus of $\mathcal{M}$ contains an
irreducible curve $Z\subset\mathbb{P}^{5}$. Then $Z\subset Y$, but
the generality of $\psi$ and the irreducibility of $Z$ and $C$
imply that $\psi(Z)=C$ and
$$
\Lambda_{j}^{i}\subset Z,
$$
but $\psi\vert_{Z}:Z\to C$ is a birational morphism. In
particular, the equality $\mathrm{deg}(Z)=j$ holds, but $Z$
contains at least $|\Lambda_{j}^{i}|$ points of
$\Sigma\subset\mathbb{P}^{4}$, which is impossible by
Proposition~\ref{proposition:factoriality-nodal-complete-intersections-nodes-in-general-position}.
\end{proof}

\begin{corollary}
\label{corollary:factoriality-nodal-complete-intersections-k-is-at-most-two}
The inequality $r\geqslant 2$ holds.
\end{corollary}

Let $\Xi_{j}^{i}$ be a base locus of the linear system of
hypersurfaces in $\mathbb{P}^{4}$ of degree $j$ that contains the
set $\Lambda_{j}^{i}$. Then $\Xi_{j}^{i}$ is a finite subset in
$\mathbb{P}^{5}$ by
Lemma~\ref{lemma:factoriality-nodal-complete-intersections-complete-intersection-III},
while we have $\Lambda_{j}^{i}\subseteq\Xi_{j}^{i}$.

\begin{lemma}
\label{lemma:factoriality-nodal-complete-intersections-complete-intersection-V}
Suppose that $\Xi_{j}^{i}\ne\varnothing$. Then the points of the
set $\Xi_{j}^{i}$ impose independent linear conditions on
hypersurfaces in $\mathbb{P}^{5}$ of degree $5(j-1)$.
\end{lemma}

\begin{proof}
The claim follows from Lemma~\ref{lemma:non-vanishing}.
\end{proof}

In particular, the points of the set $\Lambda_{j}^{i}$ impose
independent linear conditions on hypersurfaces in $\mathbb{P}^{5}$
of degree $5(j-1)$ in the case when
$\Lambda_{j}^{i}\ne\varnothing$.

\begin{lemma}
\label{lemma:factoriality-nodal-complete-intersections-complete-intersection-VI}
Suppose that $\bar{\Sigma}=\varnothing$. Then there is a
hypersurface in $\mathbb{P}^{5}$ of degree $2n+k-6$ that  contains
all points of the set $\Sigma\setminus P$ and does not contain the
point $P\in\Sigma$.
\end{lemma}

\begin{proof}
We have a disjoint union of the subsets
$$
\Sigma=\bigcup_{j=r}^{l}\bigcup_{i=1}^{c_{j}}\Lambda_{j}^{i},
$$
which implies that there is a unique set $\Lambda_{a}^{b}$ that
contains the point $P$. In particular, the point $P$ is contained
in the set $\Xi_{a}^{b}$.

It follows from
Lemma~\ref{lemma:factoriality-nodal-complete-intersections-complete-intersection-V}
that for every non-empty set $\Xi_{j}^{i}$ containing $P$ there is
a hypersurface of degree $5(j-1)$ that  passes through all points
of the set $\Xi_{j}^{i}\setminus P$ and does not contain the point
$P$. On the other hand, the construction of the set $\Xi_{j}^{i}$
implies that for every non-empty set $\Xi_{j}^{i}$ not containing
$P$ there is a hypersurface of degree $j$ that  passes through all
points of the set $\Xi_{j}^{i}$ and does not contain the point
$P$.

We have $j<5(j-1)$, because  $j\geqslant r\geqslant 2$ (see
Corollary~\ref{corollary:factoriality-nodal-complete-intersections-k-is-at-most-two}).

Thus, for every $\Xi_{j}^{i}$ containing $P$ there is hypersurface
$F_{j}^{i}\subset\mathbb{P}^{5}$ of degree $5(j-1)$ that  contains
the set $\Xi_{j}^{i}\setminus P$ and does not contain the point
$P$. Consider a hypersurface
$$
F=\bigcup_{j=r}^{l}\bigcup_{i=1}^{c_{j}}F_{j}^{i}\subset\mathbb{P}^{5}
$$
of degree $\sum_{i=r}^{l}5(i-1)c_{i}$. Then $F$ contains
$\Sigma\setminus P$ and does not contain $P$, but
$$
\mathrm{deg}\big(F\big)=\sum_{i=r}^{l}5\big(i-1\big)c_{i}<\sum_{i=r}^{l}5ic_{i}\leqslant n-1\leqslant 2n+k-6%
$$
by
Corollary~\ref{corollary:factoriality-nodal-complete-intersections-from-the-number-of-good-points},
because $n\geqslant 5$.
\end{proof}

Put $\hat{\Sigma}=\cup_{j=r}^{l}\cup_{i=1}^{c_{j}}\Lambda_{j}^{i}$
and $\check{\Sigma}=\Sigma\setminus\hat{\Sigma}$. Then
$\Sigma=\hat{\Sigma}\cup\check{\Sigma}$ and
$\psi(\check{\Sigma})=\bar{\Sigma}\subset\Pi$.

\begin{remark}
\label{remark:factoriality-nodal-complete-intersections-hypersurface-through-bad-points}
The proof of
Lemma~\ref{lemma:factoriality-nodal-complete-intersections-complete-intersection-VI}
implies the existence of a hypersurface $F\subset\mathbb{P}^{5}$
of degree $\sum_{i=r}^{l}5(i-1)c_{i}$ such that $F$ passes through
all points of the subset $\hat{\Sigma}\setminus P\subsetneq\Sigma$
and does not contain the point $P\in\Sigma$.
\end{remark}

Put $d=2n+k-6-\sum_{i=r}^{l}5(i-1)c_{i}$. Let us check that the
subset $\bar{\Sigma}\subset\Pi\cong\mathbb{P}^{2}$ and the number
$d$ satisfy all conditions of Theorem~\ref{theorem:Bese}. We may
assume that $\hat{\Sigma}\ne\varnothing$ and
$\check{\Sigma}\ne\varnothing$.

\begin{lemma}
\label{lemma:factoriality-nodal-complete-intersections-complete-intersection-VII}
The inequality $d\geqslant 5$ holds.
\end{lemma}

\begin{proof}
The claim follows from
Corollary~\ref{corollary:factoriality-nodal-complete-intersections-from-the-number-of-good-points},
because $c_{r}\geqslant 1$.
\end{proof}

\begin{lemma}
\label{lemma:factoriality-nodal-complete-intersections-complete-intersection-VIII}
The inequality $|\bar{\Sigma}|\leqslant (d^{2}+9d+10)/6$ holds.
\end{lemma}

\begin{proof}
Let us show that $6(n+k-2)(n-1-\sum_{i=r}^{l}5ic_{i})$ does not
exceed
$$
5\Big(2n+k-6-\sum_{i=r}^{l}5(i-1)c_{i}\Big)^{2}+45\Big(2n+k-6-\sum_{i=r}^{l}5(i-1)c_{i}\Big)+50,
$$
which implies $|\bar{\Sigma}|\leqslant (d^{2}+9d+10)/6$, because
$$
\big|\bar{\Sigma}\big|<{\frac{\big(n+k-2\big)}{5}}\Big(n-1-5\sum_{i=r}^{l}ic_{i}\Big)
$$
due to the
inequality~\ref{equation:factoriality-nodal-complete-intersections-number-of-good-points}.
Suppose that the inequality that we want to prove is not true, and
put $A=n-1-\sum_{i=r}^{l}5ic_{i}$ and $B=\sum_{i=r}^{l}5c_{i}$.
Then
$$
6A\big(n+k-2\big)>5\big(A+n+k-5+B\big)^{2}+45\big(A+n+k-5+B\big)+50,
$$
which is impossible, because $A>0$ by
Corollary~\ref{corollary:factoriality-nodal-complete-intersections-from-the-number-of-good-points}
and $n\geqslant 5$.
\end{proof}

\begin{lemma}
\label{lemma:factoriality-nodal-complete-intersections-complete-intersection-X}
At most $t(d+3-t)-2$ points of the set $\bar{\Sigma}$ lie on a
curve of degree $t$ for every $t\leqslant(d+3)/2$.
\end{lemma}

\begin{proof}
Suppose that $t=1$. Then
$$
t\big(d+3-t\big)-2=d=2n+k-6-\sum_{i=r}^{l}5(i-1)c_{i}\geqslant n+k-5+\sum_{i=r}^{l}5c_{i}\geqslant n+k-5+5c_{r}>n+k-2%
$$
by
Corollary~\ref{corollary:factoriality-nodal-complete-intersections-from-the-number-of-good-points},
which implies that at most $d$ points of$\bar{\Sigma}$ lie on a
line by
Proposition~\ref{proposition:factoriality-nodal-complete-intersections-nodes-in-general-position}.

Suppose that $t>1$. The points of the subset
$\bar{\Sigma}\subset\mathbb{P}^{2}$ have property $\bigstar$,
which implies that at most $(n+k-2)t$ points of $\bar{\Sigma}$ lie
on a curve of degree $t$. Therefore, we it is enough to show that
$$
t\big(d+3-t\big)-2\geqslant \big(n+k-2\big)t
$$
for all $t>1$ such that $t\leqslant (d+3)/2$ and
$t(d+3-t)-2<|\bar{\Sigma}|$.

It is easy to see that
$$
t\big(d+3-t\big)-2\geqslant t\big(n+k-2\big)\iff
n-1-\sum_{i=r}^{l}5(i-1)c_{i}>t,
$$
because $t>1$. Suppose that the inequalities
$$
n-1-\sum_{i=r}^{l}5(i-1)c_{i}\leqslant t\leqslant {\frac{d+3}{2}}
$$
and $t(d+3-t)-2<|\bar{\Sigma}|$ hold. Let us show that our
assumptions lead to a contradiction.

Put $g(x)=x(d+3-x)-2$. Then $g(x)$ is increasing for $x\leqslant
{\frac{d+3}{2}}$. Hence, we have
$$
g(t)\geqslant g\Big(n-1-\sum_{i=r}^{l}5(i-1)c_{i}\Big),
$$
which implies the inequalities
$$
{\frac{n+k-2}{5}}\Big(n-1-\sum_{i=r}^{l}5ic_{i}\Big)>|\bar{\Sigma}|>g(t)\geqslant g\Big(n-1-\sum_{i=r}^{l}5(i-1)c_{i}\Big).%
$$

Let $A=n-1-\sum_{i=r}^{l}5ic_{i}$ and $B=\sum_{i=r}^{l}5c_{i}$.
Then
$$
A{\frac{n+k-2}{5}}>g\big(A+B\big),
$$
where $A>0$ by
Corollary~\ref{corollary:factoriality-nodal-complete-intersections-from-the-number-of-good-points}.
Hence, we have
$$
0>4(n+k-2)(A+B)+5(A+B)-2\geqslant 118,
$$
which is a contradiction.
\end{proof}

There is a curve $C\subset\Pi$ of degree
$2n+k-6-\sum_{i=r}^{l}5(i-1)c_{i}$ that contains
$\bar{\Sigma}\setminus\hat{P}$ and does not contain $\hat{P}$ by
Theorem~\ref{theorem:Bese}, and there is a hypersurface $F$ of
degree $\sum_{i=r}^{l}5(i-1)c_{i}$ that contains
$\hat{\Sigma}\setminus P$ and does not contain $P$. Let $G$ be a
cone over the curve $C$ whose vertex is $\Gamma$. Then $F\cup G$
is a hypersurface of degree $2n+k-6$ that contains
$\Sigma\setminus P$ and does not contain $P$, which concludes the
proof of
Theorem~\ref{theorem:factoriality-nodal-complete-intersections-main}.

\section{Double hypersurfaces.}
\label{section:factoriality-double-hypersurfaces}

Let $\eta:X\to F$ be double cover such that $F$ is a smooth
hypersurface of degree $n\geqslant 2$, and $\eta$ is branched in a
nodal surface $S\subset F$ that is cut out on the hypersurface $F$
by a hypersurface $G\subset\mathbb{P}^{4}$ of degree $2r\geqslant
n$. In this section we prove the following result.

\begin{theorem}
\label{theorem:factoriality-double-hypersurfaces-main}%
Suppose that $|\mathrm{Sing}(X)|\leqslant(2r+n-2)r/4$. Then $X$ is
$\mathbb{Q}$-factorial.
\end{theorem}

The following result follows from Corollary~\ref{corollary:Cynk}.

\begin{proposition}
\label{proposition:factoriality-double-hypersurfaces-defect}%
The threefold $X$ is $\mathbb{Q}$-factorial if and only if the
singular points of the surface $S$ impose independent linear
conditions on the sections in
$H^{0}(\mathcal{O}_{\mathbb{P}^{4}}(3r+n-5)\vert_{F})$.
\end{proposition}

\begin{corollary}
\label{corollary:factoriality-double-hypersurfaces-defect}%
Suppose $|\mathrm{Sing}(X)|\leqslant 3r+n-4$. Then $X$ is
$\mathbb{Q}$-factorial.
\end{corollary}

Let us prove
Theorem~\ref{theorem:factoriality-double-hypersurfaces-main}.
Suppose that $|\mathrm{Sing}(X)|\leqslant(2r+n-2)r/4$. We is about
to show that the singular points of the surface
$S\subset\mathbb{P}^{4}$ impose independent linear conditions on
hypersurfaces of degree $3r-n-5$. We may assume that $r\geqslant
3$ and $n\geqslant 2$, because the claim of
Theorem~\ref{theorem:factoriality-double-hypersurfaces-main}
follows from
Corollary~\ref{corollary:factoriality-double-hypersurfaces-defect}
and \cite{Ch04t} otherwise.

\begin{lemma}
\label{lemma:factoriality-double-hypersurfaces-isolated-singularities}
There is a hypersurface $\hat{G}\subset\mathbb{P}^{4}$ of degree
$2r$ such that the surface $S$ is a complete intersection of
$\hat{G}$ and $F$, but $\mathrm{Sing}(\hat{G})\subseteq
\mathrm{Sing}(S)$.
\end{lemma}

\begin{proof}
See the proof of
Lemma~\ref{lemma:factoriality-nodal-complete-intersections-isolated-singularities}.
\end{proof}

We may assume that $\mathrm{Sing}(G)\subseteq\mathrm{Sing}(S)$.
Let $\Sigma=\mathrm{Sing}(S)$, and $P$ be an arbitrary point of
the set $\Sigma$. We must show the existence of a hypersurface of
degree $3r+n-5$ that contains $\Sigma\setminus P$ and does not
contain $P$. The proof of
Proposition~\ref{proposition:factoriality-nodal-complete-intersections-nodes-in-general-position}
implies that at most $t(2r+n-2)$ points of the set $\Sigma$ lie on
a curve in $\mathbb{P}^{4}$ of degree $t\in\mathbb{N}$.

\begin{lemma}
\label{lemma:factoriality-double-hypersurfaces-I} Suppose that
there is a plane $\Pi\subset\mathbb{P}^{4}$ such that
$\Sigma\subset\Pi$. Then there is~hyper\-sur\-face~of degree
$3r+n-5$ that contains $\Sigma\setminus P$ and does not contain
$P$.
\end{lemma}

\begin{proof}
It follows from the proof of
Lemma~\ref{lemma:factoriality-nodal-complete-intersections-complete-intersection-I}
that to conclude the proof it is enough to check that we can apply
Corollary~\ref{corollary:Bese} to $\Sigma\subset\Pi$ and the
number $d=3r+n-5\geqslant 6$.

The inequality
$$
\big|\Sigma\big|\leqslant {\frac{d^{2}+9d+16}{6}}
$$
is obvious, because  $r\geqslant 3$, $2r\geqslant n$ and
$|\Sigma|\leqslant(2r+n-2)r/4$. Therefore, we must show that at
most $t(3r+n-2-t)-2$ points of $\Sigma$ lie on a curve of degree
$t\leqslant (3r+n-2)/2$, which implies that it is enough to show
that
$$
t\big(3r+n-2-t\big)-2\geqslant t\big(2r+n-2\big)
$$
for all $t$ such that $t\leqslant (3r+n-2)/2$ and
$t(3r+n-2-t)-2<|\Sigma|$.

We may assume that $t\geqslant 2$, because $3r+n-5\geqslant
2r+n-2$. Then
$$
t\big(3r+n-2-t\big)-2\geqslant t\big(2r+n-2\big)\iff r>t.
$$

Suppose that $r\leqslant t$ for some natural $t$ such that
$$
t\leqslant {\frac{3r+n-2}{2}}
$$
and $t(3r+n-2-t)-2<|\Sigma|$. Put $g(x)=x(3r+n-2-x)-2$. Then
$g(x)$ is increasing for all $x<(3r+n-2)/2$, which implies that
$g(t)\geqslant g(r)$. Therefore, we have
$$
{\frac{(2r+n-1)r}{4}}\geqslant\big|\Sigma\big|>g(t)\geqslant g(r)=r\big(2r+n-2\big)-2,%
$$
which is impossible for $r\geqslant 3$.
\end{proof}

Let $\Pi$ and $\Gamma$ be general plane and a line in
$\mathbb{P}^{4}$ respectively, and
$\psi:\mathbb{P}^{4}\dasharrow\Pi$ be a projection from the line
$\Gamma$. Put
$\Sigma^{\prime}=\psi(\Sigma)\subset\Pi\cong\mathbb{P}^{2}$ and
$\hat{P}=\psi(P)\in\Sigma^{\prime}$.

\begin{lemma}
\label{lemma:factoriality-double-hypersurfaces-II} Suppose that at
most $t(2r+n-2)$ points of the set $\Sigma^{\prime}$ lie on a
possibly reducible curve of degree $t\in\mathbb{N}$. Then there is
hypersurface in $\mathbb{P}^{4}$ of degree $3r+n-5$ that contains
the set $\Sigma\setminus P$ and does not contain the point $P$.
\end{lemma}

\begin{proof}
The proof of Lemma~\ref{lemma:factoriality-double-hypersurfaces-I}
implies the existence of a curve $C\subset\Pi$ of degree $3r+n-5$
that contains the set $\Sigma^{\prime}\setminus\hat{P}$ and does
not contain the point $\hat{P}$. Let $Y$ be a cone in
$\mathbb{P}^{4}$ over the curve $C$ whose vertex is a line
$\Gamma$. Then $Y$ is a hypersurface in $\mathbb{P}^{4}$ of degree
$3r+n-5$ that contains the set $\Sigma\setminus P$ and does not
contain the point $P\in\Sigma$.
\end{proof}

Therefore, to conclude the proof of
Theorem~\ref{theorem:factoriality-double-hypersurfaces-main} we
may assume that the points of the set
$\Sigma^{\prime}\subset\Pi\cong\mathbb{P}^{2}$ do not satisfy the
conditions of
Lemma~\ref{lemma:factoriality-double-hypersurfaces-II}, which
implies that there is a subset $\Lambda_{k}^{1}\subset\Sigma$ such
that $|\Lambda_{k}^{1}|>k(2r+n-2)$, but the points of the set
$$
\tilde{\Lambda}_{k}^{1}=\psi(\Lambda_{k}^{1})\subset\Sigma^{\prime}\subset\Pi\cong\mathbb{P}^{2}
$$
are contained in a curve $C\subset\Pi$ of degree $k$. Moreover, we
may assume that $k$ is the minimal natural number of such
property, which implies that $C$ is irreducible and reduced.

We can iterate the construction of the subset
$\Lambda_{k}^{1}\subset\Sigma$ to get the disjoint union of
subsets
$$
\bigcup_{j=k}^{l}\bigcup_{i=1}^{c_{j}}\Lambda_{j}^{i}\subset\Sigma
$$
such that $|\Lambda_{j}^{i}|>j(2r+n-2)$, the points of the set
$\tilde{\Lambda}_{j}^{i}=\psi(\Lambda_{j}^{i})$ lie on an
irreducible curve of degree $j$, and at most $t(2r+n-2)$ points of
the subset
$$
\bar{\Sigma}=\Sigma^{\prime}\setminus\bigcup_{j=k}^{l}\bigcup_{i=1}^{c_{j}}\tilde{\Lambda}_{j}^{i}\subsetneq\Sigma^{\prime}\subset\Pi\cong\mathbb{P}^{2}
$$
are contained in a curve of degree $t$, where $c_{j}\geqslant 0$
and $c_{k}>0$. Hence, we have
\begin{equation}
\label{equation:factoriality-double-hypersurfaces-number-of-good-points}
\big|\bar{\Sigma}\big|<{\frac{\big(2r+n-2\big)r}{4}}-\sum_{i=k}^{l}c_{i}\big(2r+n-2\big)i={\frac{2r+n-2}{4}}\Big(r-\sum_{i=k}^{l}4ic_{i}\Big).
\end{equation}

\begin{corollary}
\label{corollary:factoriality-double-hypersurfaces-from-the-number-of-good-points}
The inequality $\sum_{i=k}^{l}ic_{i}<r/4$ holds.
\end{corollary}

\begin{lemma}
\label{lemma:factoriality-double-hypersurfaces-III} Let
$\mathcal{M}$ be linear system of hypersurfaces in
$\mathbb{P}^{4}$ of degree $j$ that contains all points of the set
$\Lambda_{j}^{i}$. Then the base locus of the linear system
$\mathcal{M}$ is zero-dimensional.
\end{lemma}

\begin{proof}
See the proof of
Lemma~\ref{lemma:factoriality-nodal-complete-intersections-complete-intersection-III}.
\end{proof}

\begin{corollary}
\label{corollary:factoriality-double-hypersurfaces-k-is-at-most-two}
The inequality $k\geqslant 2$ holds.
\end{corollary}

Let $\Xi_{j}^{i}$ be a base locus of the linear system of
hypersurfaces in $\mathbb{P}^{4}$ of degree $j$ that contains the
set $\Lambda_{j}^{i}$. Then $\Xi_{j}^{i}$ is a finite subset in
$\mathbb{P}^{4}$ by
Lemma~\ref{lemma:factoriality-double-hypersurfaces-III} such that
$\Lambda_{j}^{i}\subseteq\Xi_{j}^{i}$.

\begin{lemma}
\label{lemma:factoriality-double-hypersurfaces-V} The points of
the set $\Xi_{j}^{i}$ impose independent linear conditions on the
hypersurface of degree $4(j-1)$.
\end{lemma}

\begin{proof}
The required claim follows from Lemma~\ref{lemma:non-vanishing}.
\end{proof}

The points of $\Lambda_{j}^{i}$ impose independent linear
conditions on hypersurface of degree $4(j-1)$.

\begin{lemma}
\label{lemma:factoriality-double-hypersurfaces-VI} Suppose that
$\bar{\Sigma}=\varnothing$. Then there is hypersurface in
$\mathbb{P}^{4}$ degree $3r+n-5$ that  contains set
$\Sigma\setminus P$, but  does not contain the point $P\in\Sigma$.
\end{lemma}

\begin{proof}
See the proof of
Lemma~\ref{lemma:factoriality-nodal-complete-intersections-complete-intersection-VI}.
\end{proof}

Put
$\hat{\Sigma}=\cup_{j=k}^{l}\cup_{i=1}^{c_{j}}\Lambda_{j}^{i}$,
$\check{\Sigma}=\Sigma\setminus\hat{\Sigma}$ and
$d=3r+n-5-\sum_{i=k}^{l}4(i-1)c_{i}$. Then it immediately follows
from the proof of
Theorem~\ref{theorem:factoriality-nodal-complete-intersections-main}
that to conclude the proof of
Theorem~\ref{theorem:factoriality-double-hypersurfaces-main} it is
enough to check that we can apply Theorem~\ref{theorem:Bese} to
the subset $\bar{\Sigma}\subset\Pi$ and the number $d$.

\begin{lemma}
\label{lemma:factoriality-double-hypersurfaces-VII} The inequality
$d\geqslant 3$ holds.
\end{lemma}

\begin{proof}
The required claim follows from
Corollary~\ref{corollary:factoriality-double-hypersurfaces-from-the-number-of-good-points},
because $r\geqslant 3$ and $c_{k}\geqslant 1$.
\end{proof}

\begin{lemma}
\label{lemma:factoriality-double-hypersurfaces-VIII} The
inequality $|\bar{\Sigma}|\leqslant (d^{2}+9d+10)/6$ holds.
\end{lemma}

\begin{proof}
Suppose that $|\bar{\Sigma}|>(d^{2}+9d+10)/6$. Then
$$
6\big(2r+n-2\big)\Big(r-\sum_{i=k}^{l}4ic_{i}\Big)>4\big(d^{2}+9d+10\big),
$$
and putting $A=r-\sum_{i=k}^{l}4ic_{i}$ and
$B=\sum_{i=k}^{l}c_{i}$ we see that
$$
6A\big(2r+n-2\big)>4\big(2r+n-5+A+4B\big)^{2}+36\big(2r+n-5+A+4B\big)+40,
$$
where $r\geqslant 3$, but $A>0$ by
Corollary~\ref{corollary:factoriality-double-hypersurfaces-from-the-number-of-good-points},
which is a contradiction.
\end{proof}

\begin{lemma}
\label{lemma:factoriality-double-hypersurfaces-X} At most
$t(d+3-t)-2$ points of the set $\bar{\Sigma}$ lie on a possibly
reducible curve of degree $t$ for every $t\leqslant (d+3)/2$.
\end{lemma}

\begin{proof}
Let us consider the case $t=1$. Then it follows from
Corollary~\ref{corollary:factoriality-double-hypersurfaces-from-the-number-of-good-points}
that
$$
t\big(d+3-t\big)-2=d=3r+n-5-\sum_{i=k}^{l}4\big(i-1\big)c_{i}\geqslant 2r+n-5+4c_{k}\geqslant 2r+n-2.%
$$

Now we consider the case $t>1$. Then at most $(2r+n-2)t$ points of
the set $\bar{\Sigma}$ lie on a curve in $\mathbb{P}^{2}$ of
degree $t$. Therefore, to conclude the proof it is enough the show
that
$$
t\big(d+3-t\big)-2\geqslant \big(2r+n-2\big)t
$$
for every $t>1$ such that $t\leqslant (d+3)/2$ and
$t(d+3-t)-2<|\bar{\Sigma}|$. However, we have
$$
t\big(d+3-t\big)-2\geqslant t\big(2r+n-2\big)\iff r-\sum_{i=k}^{l}4\big(i-1\big)c_{i}>t,%
$$
because $t>1$. Thus, we may assume that
$$
r-\sum_{i=k}^{l}4(i-1)c_{i}\leqslant t\leqslant {\frac{d+3}{2}}
$$
and $t(d+3-t)-2<|\bar{\Sigma}|$. Let us deduce a contradiction,
which concludes the proof.

Put $g(x)=x(d+3-x)-2$. Then $g(x)$ is increasing for
$x\leqslant(d+3)/2$. Hence, we have
$$
g(t)\geqslant g\Big(r-\sum_{i=k}^{l}4(i-1)c_{i}\Big).
$$

Put $A=r-\sum_{i=k}^{l}4ic_{i}$ and $B=\sum_{i=k}^{l}c_{i}$. Then
$$
A{\frac{2r+n-2}{4}}>g\big(A+4B\big)=\big(A+4B\big)\big(2r+n-5\big)-2,
$$
which is impossible, because  $A>0$ by
Corollary~\ref{corollary:factoriality-double-hypersurfaces-from-the-number-of-good-points}.
\end{proof}

Thus, we proved that we can apply Theorem~\ref{theorem:Bese} to
the subset $\bar{\Sigma}\subset\Pi$ and the natural number $d$,
which concludes the proof of
Theorems~\ref{theorem:factoriality-double-hypersurfaces-main} and
\ref{theorem:main}.

\bigskip

\end{document}